\documentclass[12pt]{amsart}
\usepackage{amsmath, amsthm, graphicx, amssymb,verbatim,pst-all,color,txfonts}
\usepackage[all]{xy}
\usepackage[margin=0.94in]{geometry}
\usepackage{mathptmx}
\ifx\JPicScale\undefined\def\JPicScale{1.0}\fi
\unitlength \JPicScale mm
{\theoremstyle{plain}
   \newtheorem*{maintheorem}{Theorem}

}
{\theoremstyle{definition}
   \newtheorem*{ack}{Acknowledgements}

}

\newcommand{\M}{\overline{M}}

\newcommand{\PP}{\mathbb{P}}

\numberwithin{theorem}{section}

\begin{document}
\title{The dual complex of $\M_{0,n}$ via phylogenetics}
\author{Noah Giansiracusa}
\maketitle

%%%%

\begin{abstract}
The moduli space $\M_{0,n}$ of stable rational $n$-pointed curves has divisorial boundary with simple normal crossings.  In this brief note I observe that the dual complex is a flag complex; that is, a collection of boundary divisors has nonempty intersection if and only if the pairwise intersections are nonempty.  Rather than proving this directly, I translate the statement to a setting in phylogenetics where it is widely used and multiple explicit proofs have been written.  It appears this result is known by experts but lacks a detailed reference in the literature, except recently for $n=7$.  
\end{abstract} 

%%%%

\renewcommand{\thefootnote}{\fnsymbol{footnote}} 
\footnotetext{MSC-class: 14H10, 05C05}     
\renewcommand{\thefootnote}{\arabic{footnote}} 

\section{Introduction}

The moduli space $\M_{0,n}$ of stable rational $n$-pointed curves has featured in many areas of mathematics.  It was first envisioned by Grothendieck \cite{SGA7-I}, then constructed as an iterated blow-up by Knudsen \cite{Knudsen2}, then given slicker iterated blow-up constructions by Kapranov \cite{Kapranov-chow,Kapranov-veronese} and Keel \cite{Keel-thesis}.  It has been the subject of interest in birational geometry \cite{Keel-McKernan,Castravet-Tevelev-mds}, arithmetic \cite{Goncharov-Manin}, topology \cite{realcohomology}, representation theory \cite{M0n-Gaudin}, enumerative geometry \cite{QuantumM0n}, derived categories \cite{Manin-Smirnov-derived}, etc.  Its scope even reaches beyond pure mathematics, as the tree structures arising in $\M_{0,n}$ play an important role in phylogenetics \cite{Billera-Holmes-Vogtmann,Kapranov-assoc,Speyer-Sturmfels-tropical-grassmannian,Levy-Pachter}.  In this brief note I exploit this latter connection to provide a succinct proof, via reference to a combinatorial result in the phylogenetics literature, of a folklore result about $\M_{0,n}$: a collection of boundary divisors has nonempty intersection if and only if all pairs in the collection have nonempty intersection.  This can be stated formally as follows (see \cite{Kollar-prescribed} for a nice discussion of dual complexes and their applications and \cite{QuantumM0n} for an elementary introduction to $\M_{0,n}$, and recall that a simplicial complex is called a flag complex if it is the maximal simplicial complex on its 1-skeleton):

\begin{maintheorem}
The dual complex of the divisorial boundary in $\M_{0,n}$ is a flag complex.
\end{maintheorem}

The really is a special property; for instance, the three coordinate lines in $\PP^2$ intersect in pairs whereas their triple intersection is empty.  This fact about $\M_{0,n}$ seems to be known by some experts (cf., \cite[p.871]{Gibney-Maclagan-quotients}, \cite[Lemma 7.2]{Hacking-Keel-Tevelev}, \cite{Speyer-Sturmfels-tropical-grassmannian}) and was recently proven explicitly in the first non-trivial case, $n=7$ \cite[Lemma 3.2 and Proposition 3.4]{Luca}.  As I show below, by considering the well-known stratification of $\M_{0,n}$ by boundary strata corresponding to dual graph trees, this intersection-theoretic property of $\M_{0,n}$ is precisely the combinatorial ``Splits-Equivalence Theorem'' from the 70s due to Buneman \cite{Buneman71,Buneman}, which in turn is the ``Four Gamete Condition'' in phylogenetics \cite{Gamete1,Gamete2,Gamete3}.  Buneman's result was given an alternate proof in the 80s \cite{Bandelt-Dress} which is nicely exposed in Semple and Steel's book on mathematical phylogenetics \cite{Semple-Steel}.  

\begin{ack}
I thank Luca Schaffler for bringing to my attention this question about $\M_{0,n}$ and I thank Luca, Jenia Tevelev, and James McKernan for helpful discussions.  I also thank the Simons Foundation for financial support.
\end{ack}

%\vspace{-0.1in}

\section{Proof via reinterpretation}

A \emph{phylogenetic tree on $n$ taxa} is a graph that is a tree with $n$ leaves, labelled by the integers $1,\ldots,n$, such that there are no vertices of degree two.  Each point of $\M_{0,n}$ determines a phylogenetic tree \cite[\S4]{Kapranov-assoc}.  Indeed, the points of $\M_{0,n}$ are in bijection with stable rational $n$-pointed curves $(C,p_1,\ldots,p_n)$, which means $C$ is a nodal union of $\PP^1$s such that each irreducible component contains at least three \emph{special} points (nodes or marked points $p_i$) and such that the dual graph is a tree.  The \emph{dual graph} here refers to the graph with a vertex for each irreducible component $C_i \subseteq C$ and each marked point $p_i$, and an edge joining two vertices if they correspond to irreducible components satisfying $C_i \cap C_j \ne \varnothing$ or to a component and a marked point satisfying $p_j \in C_i$.  Such an object is necessarily a phylogenetic tree.  These dual graphs induce a widely-used stratification of $\M_{0,n}$, where a \emph{stratum} is the closure of the locus of points with fixed dual graph.  The strata in $\M_{0,n}$ are in bijection with the phylogenetic trees on $n$ taxa, where the bijection sends a stratum to the dual graph of any point in its relative interior (which we shall refer to simply as the dual graph of the stratum).

In the phylogenetics literature, a bipartition $I\sqcup J= \{1,\ldots,n\}$ is referred to as a \emph{split}.  Each \emph{internal} edge (that is, edge between non-leaf vertices) of a phylogenetic tree induces a split by considering the sets of integer labels on either side of the edge (which is well-defined due to the tree property of the graph).  The set of all such splits of a phylogenetic tree $T$ is denoted $\Sigma(T)$.  A pair of splits $I\sqcup J$, $K\sqcup L$ is said to be \emph{compatible} if at least one of the following intersections is empty (see \cite[Chapter 3]{Semple-Steel}): $I\cap K,~I\cap L, ~J\cap K, ~J\cap L$.  The famous Buneman Splits-Equivalence Theorem alluded to in the introduction is the following:

\begin{maintheorem}[\cite{Semple-Steel}, Theorem 3.1.4]
A collection of splits of $\{1,\ldots,n\}$ is of the form $\Sigma(T)$ for a phylogenetic tree $T$ on $n$ taxa if and only if the splits are pairwise compatible.
\end{maintheorem}

The irreducible boundary divisors in $\M_{0,n}$ are the codimension one strata given by dual graphs with a single internal edge; these are in bijection with splits $I\sqcup J= \{1,\ldots,n\}$ such that $I$ and $J$ have cardinality at least two.  Let us denote by $D_{I,J}$ the divisor corresponding to the split $I\sqcup J= \{1,\ldots,n\}$, so that the generic point of $D_{I,J}$ corresponds to a two-component curve with a node separating the points indexed by $I$ from those indexed by $J$.  A pair of boundary divisors has nonempty intersection, $D_{I,J}\cap D_{K,L} \ne \varnothing$, if and only if at least one of the following containments holds \cite[pp.550--551]{Keel-thesis}: $I \subset K,~ I\subset L,~ K\subset I,~ L\subset I.$  Since $J = I^c$ and $L = K^c$, a pair of boundary divisors has nonempty intersection if and only if the corresponding splits are compatible.

Each stratum in $\M_{0,n}$ is the complete intersection of the boundary divisors containing it, and conversely any intersection of boundary divisors is a single stratum if it is nonempty.  If a stratum has dual graph $T$, then the boundary divisors containing it are precisely the ones corresponding to splits in $\Sigma(T)$.  Putting this together, we see that a collection of boundary divisors has nonempty intersection if and only if the corresponding set of splits is of the form $\Sigma(T)$ for some phylogenetic tree $T$, and by Buneman's result this is equivalent to the splits being pairwise compatible, which in turn is equivalent to the boundary divisors having pairwise nonempty intersection.  \hfill$\Box$

%%%%

\large{\bibliographystyle{amsalpha}}
\bibliography{bib}

\end{document}